\documentclass[12pt,twoside]{article}

\input{prepictex.tex}
\input{pictex.tex}
\input{postpictex.tex}

\usepackage[utf8]{inputenc}
\usepackage[english]{babel}
\usepackage{hyperref}
\usepackage{amsthm,thmtools}
\usepackage{amsfonts}
\usepackage{amssymb}
\usepackage{amsmath}
\usepackage{float}
\usepackage{pgf}
\usepackage{tikz}
\usetikzlibrary{calc}
\usetikzlibrary{arrows,decorations.pathmorphing,backgrounds,positioning,fit,petri,patterns,shapes,automata}
\usepackage{xcolor}
\usepackage{latexsym}
\newcounter{enunciato}[section]
\newtheorem{ittheorem}{Theorem}

\newenvironment{theorem}{\addtocounter{enunciato}{1}
\begin{ittheorem}}{\end{ittheorem}}
\newenvironment{corollary}{\addtocounter{enunciato}{1}
\begin{itcorollary}}{\end{itcorollary}}
\newtheorem{itlemma}{Lemma}

\newtheorem{itcorollary}{Corollary}

\usepackage[figurename=Fig]{caption}
\usepackage[labelfont=bf]{caption}
\usepackage[labelsep=period]{caption}
\usepackage{caption}
\usepackage{subcaption}
\usepackage{graphicx}
\usepackage{geometry}
\usepackage{subfiles}
\usepackage{float}
\usepackage{mathtools}
\usepackage{pgfplots}
\pgfplotsset{/pgf/number format/use comma,compat=newest}
\usepackage{url}
\usepackage{authblk}

\title{A ${\mathbb N}$atural Avenue}
\author[1]{Roberto Conti}
\author[2]{Pierluigi Contucci}
\affil[1]{Department of Basic and Applied Sciences for Engineering, Sapienza University of Rome}
\affil[2]{Department of Mathematics, University of Bologna}
\date{\today}

\begin{document}

\maketitle

\begin{abstract}
We consider an infinite sequence of rooted trees 
naturally 
emerging in a number-theoretical context. We advance some ideas on its structure by discussing some elementary properties. Some of those properties are shown to be related to classical results or conjectures in number theory.
\end{abstract}

\section{Introduction}
The fundamental theorem of arithmetic states that every natural number $n$ larger than 1 can be represented 
as a product of distinct primes raised to natural powers: 
\begin{equation}\label{pnd}
n = {p_1}^{n_1} \cdots {p_k}^{n_k} .
\end{equation}
Such representation is not unique, but it becomes so if we further require that the distinct primes appear in increasing order or, otherwise said, it's unique up to a permutation of the factors. Note that the integer $k$ as above, i.e. the number of distinct prime factors of $n$, is an arithmetic function denoted $\omega(n)$ in the number theoretical literature.

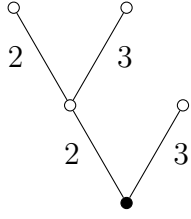
\begin{figure}
    \centering
    \begin{tikzpicture}
\tikzstyle{hollow node}=[circle,draw,inner sep=1.5]
\tikzstyle{solid node}=[circle,draw,inner sep=1.5,fill=black]
\node(0)[solid node]{}
child[grow=60]{node[hollow node]{} edge from parent node[xshift=10]{$3$}}
child[grow=120]{node[hollow node]{}  
child[grow=60]{node[hollow node]{} edge from parent node[xshift=10]{$3$}}
child[grow=120]{node[hollow node]{} edge from parent node[xshift=-10]{$2$}}
edge from parent node[xshift=-10]{$2$}
};
\end{tikzpicture}
    \caption{Labeled tree representation of the integer $192=2^{2\cdot 3}\cdot 3$}
    \label{fig:my_label}
\end{figure}

Looking for a pure representation of each integer {\it only} in terms of prime numbers (the {\it building blocks} of arithmetic) we proceed by iterating such decomposition to the successive natural powers $n_i$ in (\ref{pnd}). Some examples are $6=2\cdot 3$, 
$12=2^2\cdot 3$, $18=2\cdot 3^2$, $486=2\cdot 3^5$ and
\begin{equation}\label{154}
1549681956=2^2\cdot 3^{2\cdot 3^2} \, .
\end{equation}
The previous representation identifies, for each integer $n$, a graph theoretical structure which is easily seen to be a planar rooted tree: starting with a given marked vertex called the root (somehow representing the number $1$) we draw $k=\omega(n)$ edges emanating from the root, each labeled with one of the primes $p_i$ as above. We draw the tree starting from the root at the bottom and proceed to the top like in real trees. We then apply the same procedure to the endpoints of each of these edges, say the one labeled by $p_i$, to which we attach new $k_i=\omega(n_i)$ edges corresponding to the distinct primes appearing in the decomposition of $n_i$, and so on.
By iterating this procedure, as many times as needed, we end up with a finite rooted planar tree with edges labeled by suitable primes and uniquely associated to each integer. By construction we see that not all labelings are allowed, but they are subject to the only restriction that each branching 
involves only distinct primes. The labeled tree corresponding to the integer $192$ is shown in Fig \ref{fig:my_label}. The representation of an integer with powers of primes can also be also found in \cite{DG} where is called prime tower factorization (see also \cite{VVI1}).

\medskip
Let us consider the set of all finite rooted planar trees ${\cal T}$ (see \cite{St}). Then the previous construction identifies a canonical map $t$
\begin{equation}
{\mathbb N} \to {\cal T} \; ,
\end{equation}
by simply ignoring the labeling primes. 

We observe moreover that by identifying trees with the same topological structure, no matter how they are embedded in the plane, and calling such set $[{\cal T}]$, we also get a quotient map $t^\#$
\begin{equation}
{\mathbb N} \to [{\cal T}] \; .
\end{equation}

\begin{figure}
     
     \begin{subfigure}[b]{0.24\textwidth}
     \centering
\begin{tikzpicture}
\tikzstyle{hollow node}=[circle,draw,inner sep=1.5]
\tikzstyle{solid node}=[circle,draw,inner sep=1.5,fill=black]
\node(0)[solid node]{}
child[grow=90]{node[hollow node]{}};
\end{tikzpicture}
\caption{$t(13)=(())$}
    \label{13}
     \end{subfigure}
     \hfill
     \begin{subfigure}[b]{0.24\textwidth}
     \centering
\begin{tikzpicture}
\tikzstyle{hollow node}=[circle,draw,inner sep=1.5]
\tikzstyle{solid node}=[circle,draw,inner sep=1.5,fill=black]
\node(0)[solid node]{}
child[grow=60]{node[hollow node]{}}
child[grow=120]{node[hollow node]{}};
\end{tikzpicture}
\caption{$t(14)=(()())$}
    \label{14}
     \end{subfigure}
     \hfill
     \begin{subfigure}[b]{0.24\textwidth}
     \centering
\begin{tikzpicture}
\tikzstyle{hollow node}=[circle,draw,inner sep=1.5]
\tikzstyle{solid node}=[circle,draw,inner sep=1.5,fill=black]
\node(0)[solid node]{}
child[grow=60]{node[hollow node]{}}
child[grow=120]{node[hollow node]{}};
\end{tikzpicture}
\caption{$t(15)=(()())$}
    \label{15}
     \end{subfigure}
     \hfill
     \begin{subfigure}[b]{0.24\textwidth}
     \centering
     \vspace{0.5cm}
\begin{tikzpicture}
\tikzstyle{hollow node}=[circle,draw,inner sep=1.5]
\tikzstyle{solid node}=[circle,draw,inner sep=1.5,fill=black]
\node(0)[solid node]{}
child[grow=90]{node[hollow node]{}child[grow=90]{node[hollow node]{}child[grow=90]{node[hollow node]{}}}};
\end{tikzpicture}
\caption{$t(16)=(((())))$}
    \label{16}
     \end{subfigure}

      \begin{subfigure}[b]{0.24\textwidth}
     \centering
\begin{tikzpicture}
\tikzstyle{hollow node}=[circle,draw,inner sep=1.5]
\tikzstyle{solid node}=[circle,draw,inner sep=1.5,fill=black]
\node(0)[solid node]{}
child[grow=90]{node[hollow node]{}};
\end{tikzpicture}
\caption{$t(17)=(())$}
    \label{17}
     \end{subfigure}
     \hfill
     \begin{subfigure}[b]{0.24\textwidth}
     \centering
\vskip 1truecm
\begin{tikzpicture}
\tikzstyle{hollow node}=[circle,draw,inner sep=1.5]
\tikzstyle{solid node}=[circle,draw,inner sep=1.5,fill=black]
\node(0)[solid node]{}
child[grow=60]{node[hollow node]{}child[grow=60]{node[hollow node]{}}}
child[grow=120]{node[hollow node]{}};
\end{tikzpicture}
\caption{$t(18)=(()(()))$}
    \label{18}
     \end{subfigure}
     \hfill
     \begin{subfigure}[b]{0.24\textwidth}
     \centering
\begin{tikzpicture}
\tikzstyle{hollow node}=[circle,draw,inner sep=1.5]
\tikzstyle{solid node}=[circle,draw,inner sep=1.5,fill=black]
\node(0)[solid node]{}
child[grow=90]{node[hollow node]{}};
\end{tikzpicture}
\caption{$t(19)=(())$}
    \label{19}
     \end{subfigure}
     \hfill
     \begin{subfigure}[b]{0.24\textwidth}
     \centering
\begin{tikzpicture}
\tikzstyle{hollow node}=[circle,draw,inner sep=1.5]
\tikzstyle{solid node}=[circle,draw,inner sep=1.5,fill=black]
\node(0)[solid node]{}
child[grow=60]{node[hollow node]{}}
child[grow=120]{node[hollow node]{}child[grow=120]{node[hollow node]{}}};
\end{tikzpicture}
\caption{$t(20)=((())())$}
    \label{20}
     \end{subfigure}

      \begin{subfigure}[b]{0.24\textwidth}
     \centering
\begin{tikzpicture}
\tikzstyle{hollow node}=[circle,draw,inner sep=1.5]
\tikzstyle{solid node}=[circle,draw,inner sep=1.5,fill=black]
\node(0)[solid node]{}
child[grow=60]{node[hollow node]{}}
child[grow=90]{node[hollow node]{}}
child[grow=120]{node[hollow node]{}};
\end{tikzpicture}
\caption{$t(30)=(()()())$}
    \label{30}
     \end{subfigure}
     \hfill
     \begin{subfigure}[b]{0.24\textwidth}
     \centering
\vskip 1truecm
\begin{tikzpicture}
\tikzstyle{hollow node}=[circle,draw,inner sep=1.5]
\tikzstyle{solid node}=[circle,draw,inner sep=1.5,fill=black]
\node(0)[solid node]{}
child[grow=60]{node[hollow node]{}}
child[grow=120]{node[hollow node]{}child[grow=120]{node[hollow node]{}child[grow=120]{node[hollow node]{}}}};
\end{tikzpicture}
\caption{$t(48)=(((()))())$}
    \label{48}
     \end{subfigure}
     \hfill
     \begin{subfigure}[b]{0.24\textwidth}
     \centering
\begin{tikzpicture}
\tikzstyle{hollow node}=[circle,draw,inner sep=1.5]
\tikzstyle{solid node}=[circle,draw,inner sep=1.5,fill=black]
\node(0)[solid node]{}
child[grow=120]{node[hollow node]{}child[grow=120]{node[hollow node]{}}}
child[grow=90]{node[hollow node]{}}
child[grow=60]{node[hollow node]{}};
\end{tikzpicture}
\caption{$t(60)=((())()())$}
    \label{60}
     \end{subfigure}
     \hfill
     \begin{subfigure}[b]{0.24\textwidth}
     \centering
\begin{tikzpicture}
\tikzstyle{hollow node}=[circle,draw,inner sep=1.5]
\tikzstyle{solid node}=[circle,draw,inner sep=1.5,fill=black]
\node(0)[solid node]{}
child[grow=90]{node[hollow node]{}child[grow=60]{node[hollow node]{}}
child[grow=120]{node[hollow node]{}}};
\end{tikzpicture}
\caption{$t(64)=((()()))$}
    \label{64}
     \end{subfigure}
 \caption{More examples of the map $t(n)$.}
 \label{moretiofn}
\end{figure}
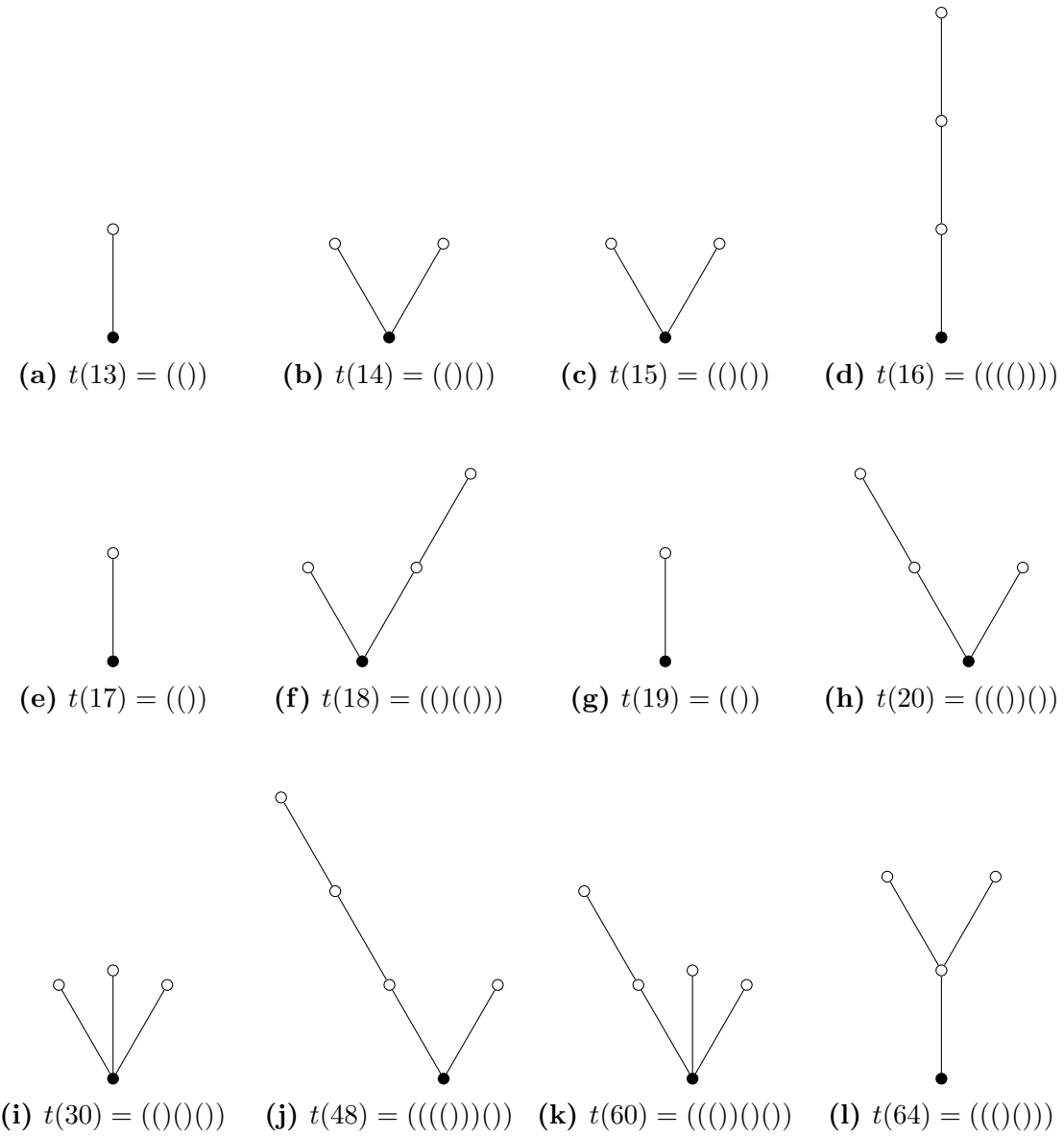

\begin{figure}
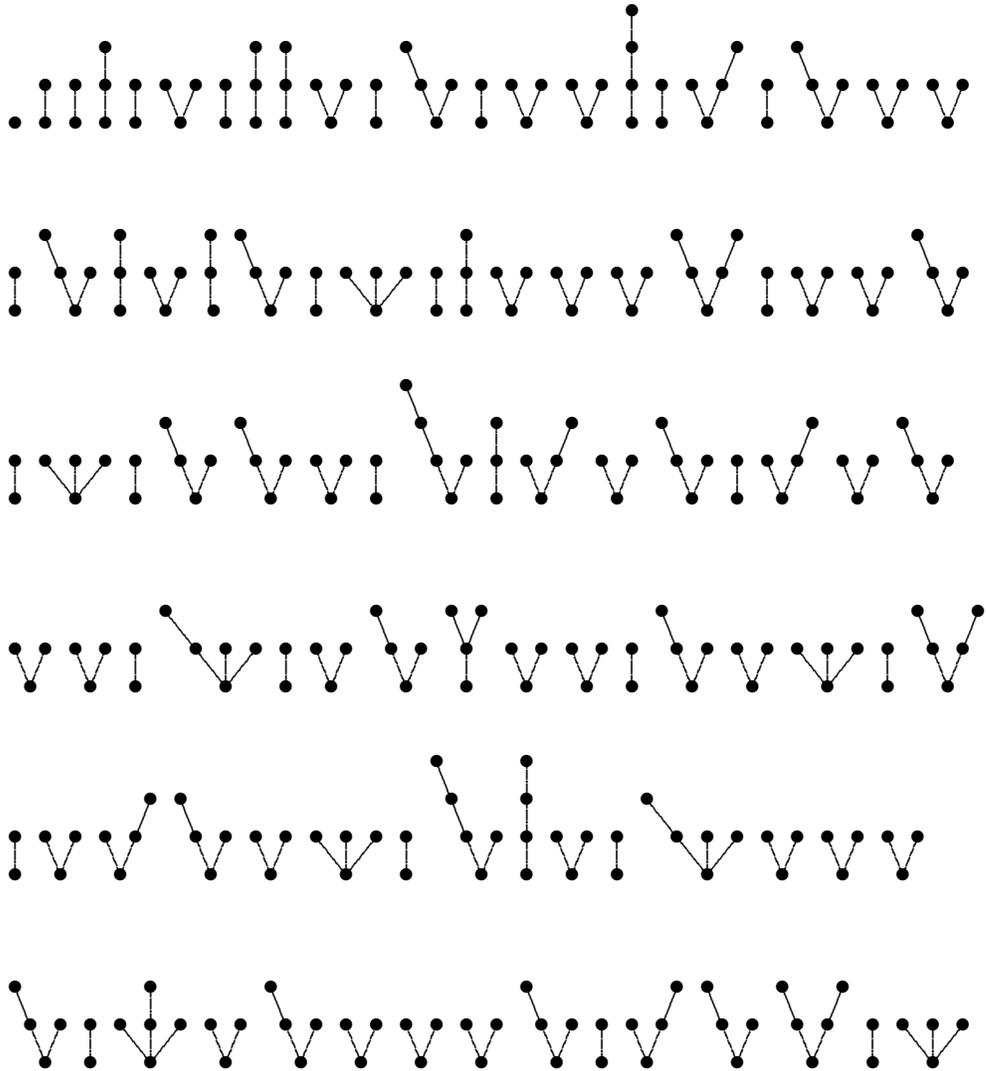

\[ \beginpicture
\setcoordinatesystem units <0.4cm,0.5cm>
\setplotarea x from 0 to 1000, y from 0 to 3

%\put {$\bigstar$} at 0 -1

\setlinear

% 1
\put {$\bullet$} at 0 0

% 2
\put {$\bullet$} at 1 0
\put {$\bullet$} at 1 1
\plot 1 0 1 1 /

% 3
\put {$\bullet$} at 2 0
\put {$\bullet$} at 2 1
\plot 2 0 2 1 /

% 4 = 2^2
\put {$\bullet$} at 3 0
\put {$\bullet$} at 3 1
\put {$\bullet$} at 3 2
\plot 3 0 3 1 /
\plot 3 1 3 2 /

% 5
\put {$\bullet$} at 4 0
\put {$\bullet$} at 4 1
\plot 4 0 4 1 /

% 6 = 2*3
\put {$\bullet$} at 5.5 0
 \put {$\bullet$} at 5 1
 \put {$\bullet$} at 6 1
\plot 5.5 0 5 1 /
\plot 5.5 0 6 1 /

% 7
\put {$\bullet$} at 7 0
\put {$\bullet$} at 7 1
\plot 7 0 7 1 /

% 8 = 2^3
\put {$\bullet$} at 8 0
\put {$\bullet$} at 8 1
\put {$\bullet$} at 8 2
\plot 8 0 8 1 /
\plot 8 1 8 2 /

% 9 = 3^2
\put {$\bullet$} at 9 0
\put {$\bullet$} at 9 1
\put {$\bullet$} at 9 2
\plot 9 0 9 1 /
\plot 9 1 9 2 /

% 10 = 2*5
\put {$\bullet$} at 10.5 0
 \put {$\bullet$} at 10 1
 \put {$\bullet$} at 11 1
\plot 10.5 0 10 1 /
\plot 10.5 0 11 1 /

% 11
\put {$\bullet$} at 12 0
\put {$\bullet$} at 12 1
\plot 12 0 12 1 /

% 12 = 2^2 * 3
\put {$\bullet$} at 14 0
 \put {$\bullet$} at 13.5 1
 \put {$\bullet$} at 13 2
 \put {$\bullet$} at 14.5 1
\plot 14 0 13.5 1 /
\plot 13.5 1 13 2 /
\plot 14 0 14.5 1 /

% 13
\put {$\bullet$} at 15.5 0
\put {$\bullet$} at 15.5 1
\plot 15.5 0 15.5 1 /

% 14 = 2*7
\put {$\bullet$} at 17 0
 \put {$\bullet$} at 16.5 1
 \put {$\bullet$} at 17.5 1
\plot 17 0 16.5 1 /
\plot 17 0 17.5 1 /

% 15 = 3*5
\put {$\bullet$} at 19 0
 \put {$\bullet$} at 18.5 1
 \put {$\bullet$} at 19.5 1
\plot 19 0 18.5 1 /
\plot 19 0 19.5 1 /

% 16 = 2^{2^2}
\put {$\bullet$} at 20.5 0
\put {$\bullet$} at 20.5 1
\put {$\bullet$} at 20.5 2
\put {$\bullet$} at 20.5 3
\plot 20.5 0 20.5 1 /
\plot 20.5 1 20.5 2 /
\plot 20.5 2 20.5 3 /

% 17
\put {$\bullet$} at 21.5 0
\put {$\bullet$} at 21.5 1
\plot 21.5 0 21.5 1 /

% 18 = 2 * 3^2
\put {$\bullet$} at 23 0
 \put {$\bullet$} at 22.5 1
 \put {$\bullet$} at 23.5 1
 \put {$\bullet$} at 24 2
\plot 23 0 22.5 1 /
\plot 23 0 23.5 1 /
\plot 23.5 1 24 2 /

% 19
\put {$\bullet$} at 25 0
\put {$\bullet$} at 25 1
\plot 25 0 25 1 /

% 20 = 2^2 * 5
\put {$\bullet$} at 27 0
 \put {$\bullet$} at 26.5 1
 \put {$\bullet$} at 26 2
 \put {$\bullet$} at 27.5 1
\plot 27 0 26.5 1 /
\plot 27 0 27.5 1 /
\plot 26.5 1 26 2 /

% 21 = 3*7
\put {$\bullet$} at 29 0
 \put {$\bullet$} at 28.5 1
 \put {$\bullet$} at 29.5 1
\plot 29 0 28.5 1 /
\plot 29 0 29.5 1 /

% 22 = 2 * 11
\put {$\bullet$} at 31 0
 \put {$\bullet$} at 30.5 1
 \put {$\bullet$} at 31.5 1
\plot 31 0 30.5 1 /
\plot 31 0 31.5 1 /

% 23
\put {$\bullet$} at 0 -5
 \put {$\bullet$} at 0 -4
\plot 0 -5 0 -4 /

% 24 = 2^3 * 3
\put {$\bullet$} at 2 -5
 \put {$\bullet$} at 1.5 -4
 \put {$\bullet$} at 1 -3
 \put {$\bullet$} at 2.5 -4
\plot 2 -5 1.5 -4 /
\plot 1.5 -4 1 -3 /
\plot 2 -5 2.5 -4 /

% 25 = 5^2
\put {$\bullet$} at 3.5 -5
 \put {$\bullet$} at 3.5 -4
 \put {$\bullet$} at 3.5 -3
\plot 3.5 -5 3.5 -4 /
\plot 3.5 -4 3.5 -3 /

% 26 = 2 * 13
\put {$\bullet$} at 5 -5
 \put {$\bullet$} at 4.5 -4
 \put {$\bullet$} at 5.5 -4
\plot 5 -5 4.5 -4 /
\plot 5 -5 5.5 -4 /

% 27 = 3^3
\put {$\bullet$} at 6.6 -5
 \put {$\bullet$} at 6.5 -4
 \put {$\bullet$} at 6.5 -3
\plot 6.5 -5 6.5 -4 /
\plot 6.5 -4 6.5 -3 /

% 28 = 2^2 * 7
\put {$\bullet$} at 8.5 -5
 \put {$\bullet$} at 8 -4
 \put {$\bullet$} at 7.5 -3
 \put {$\bullet$} at 9 -4
\plot 8.5 -5 8 -4 /
\plot 8 -4 7.5 -3 /
\plot 8.5 -5 9 -4 /

% 29
\put {$\bullet$} at 10 -5
\put {$\bullet$} at 10 -4
\plot 10 -5 10 -4 /

% 30 = 2*3*5
\put {$\bullet$} at 12 -5
\put {$\bullet$} at 11 -4
\put {$\bullet$} at 12 -4
\put {$\bullet$} at 13 -4
\plot 12 -5 11 -4 /
\plot 12 -5 12 -4 /
\plot 12 -5 13 -4 /

% 31
\put {$\bullet$} at 14 -5
 \put {$\bullet$} at 14 -4
\plot 14 -5 14 -4 /

% 32 = 2^5
\put {$\bullet$} at 15 -5
 \put {$\bullet$} at 15 -4
  \put {$\bullet$} at 15 -3
\plot 15 -5 15 -4 /
\plot 15 -4 15 -3 /

% 33 = 3*11
\put {$\bullet$} at 16.5 -5
 \put {$\bullet$} at 16 -4
 \put {$\bullet$} at 17 -4
\plot 16.5 -5 16 -4 /
\plot 16.5 -5 17 -4 /

% 34 = 2*17
\put {$\bullet$} at 18.5 -5
 \put {$\bullet$} at 18 -4
 \put {$\bullet$} at 19 -4
\plot 18.5 -5 18 -4 /
\plot 18.5 -5 19 -4 /

% 35 = 5*7
\put {$\bullet$} at 20.5 -5
 \put {$\bullet$} at 20 -4
 \put {$\bullet$} at 21 -4
\plot 20.5 -5 20 -4 /
\plot 20.5 -5 21 -4 /

% 36 = 2^2 * 3^2
\put {$\bullet$} at 23 -5
 \put {$\bullet$} at 22.5 -4
  \put {$\bullet$} at 22 -3
 \put {$\bullet$} at 23.5 -4
 \put {$\bullet$} at 24 -3
\plot 23 -5 22.5 -4 /
\plot 22.5 -4 22 -3 /
\plot 23 -5 23.5 -4 /
\plot 23.5 -4 24 -3 /

% 37
\put {$\bullet$} at 25 -5
 \put {$\bullet$} at 25 -4
\plot 25 -5 25 -4 /

% 38 = 2*19
\put {$\bullet$} at 26.5 -5
 \put {$\bullet$} at 26 -4
 \put {$\bullet$} at 27 -4
\plot 26.5 -5 26 -4 /
\plot 26.5 -5 27 -4 /

% 39 = 3*13
\put {$\bullet$} at 28.5 -5
 \put {$\bullet$} at 28 -4
 \put {$\bullet$} at 29 -4
\plot 28.5 -5 28 -4 /
\plot 28.5 -5 29 -4 /

% 40 = 2^3 * 5
\put {$\bullet$} at 31 -5
 \put {$\bullet$} at 30.5 -4
 \put {$\bullet$} at 30 -3
 \put {$\bullet$} at 31.5 -4
\plot 31 -5 30.5 -4 /
\plot 30.5 -4 30 -3 /
\plot 31 -5 31.5 -4 /

% 41
\put {$\bullet$} at 0 -10
 \put {$\bullet$} at 0 -9
\plot 0 -10 0 -9 /

% 42 = 2*3*7
\put {$\bullet$} at 2 -10
\put {$\bullet$} at 1 -9
\put {$\bullet$} at 2 -9
\put {$\bullet$} at 3 -9
\plot 2 -10 1 -9 /
\plot 2 -10 2 -9 /
\plot 2 -10 3 -9 /

% 43
\put {$\bullet$} at 4 -10
 \put {$\bullet$} at 4 -9
\plot 4 -10 4 -9 /

% 44 = 2^2 * 11
\put {$\bullet$} at 6 -10
 \put {$\bullet$} at 5.5 -9
 \put {$\bullet$} at 5 -8
 \put {$\bullet$} at 6.5 -9
\plot 6 -10 5.5 -9 /
\plot 5.5 -9 5 -8 /
\plot 6 -10 6.5 -9 /

% 45 = 3^2 * 5
\put {$\bullet$} at 8.5 -10
 \put {$\bullet$} at 8 -9
 \put {$\bullet$} at 7.5 -8
 \put {$\bullet$} at 9 -9
\plot 8.5 -10 8 -9 /
\plot 8 -9 7.5 -8 /
\plot 8.5 -10 9 -9 /

% 46 = 2*23
\put {$\bullet$} at 10.5 -10
 \put {$\bullet$} at 10 -9
 \put {$\bullet$} at 11 -9
\plot 10.5 -10 10 -9 /
\plot 10.5 -10 11 -9 /

% 47
\put {$\bullet$} at 12 -10
 \put {$\bullet$} at 12 -9
\plot 12 -10 12 -9 /

% 48 = 2^{2^2} * 3
\put {$\bullet$} at 14.5 -10
\put {$\bullet$} at 14 -9
\put {$\bullet$} at 13.5 -8
\put {$\bullet$} at 13 -7
\put {$\bullet$} at 15 -9
\plot 14.5 -10 14 -9 /
\plot 14 -9 13.5 -8 /
\plot 13.5 -8 13 -7 /
\plot 14.5 -10 15 -9 /

% 49 = 7^2
\put {$\bullet$} at 16 -10
\put {$\bullet$} at 16 -9
\put {$\bullet$} at 16 -8
\plot 16 -10 16 -9 /
\plot 16 -9 16 -8 /

% 50 = 2 * 5^2
\put {$\bullet$} at 17.5 -10
 \put {$\bullet$} at 17 -9
 \put {$\bullet$} at 18 -9
 \put {$\bullet$} at 18.5 -8
\plot 17.5 -10 17 -9 /
\plot 17.5 -10 18 -9 /
\plot 18 -9 18.5 -8 /

% 51 = 3*17
\put {$\bullet$} at 20 -10
 \put {$\bullet$} at 19.5 -9
 \put {$\bullet$} at 20.5 -9
\plot 20 -10 19.5 -9 /
\plot 20 -10 20.5 -9 /

% 52 = 2^2 * 13
\put {$\bullet$} at 22.5 -10
 \put {$\bullet$} at 22 -9
 \put {$\bullet$} at 21.5 -8
 \put {$\bullet$} at 23 -9
\plot 22.5 -10 22 -9 /
\plot 22 -9 21.5 -8 /
\plot 22.5 -10 23 -9 /

% 53
\put {$\bullet$} at 24 -10
 \put {$\bullet$} at 24 -9
\plot 24 -10 24 -9 /

% 54 = 2 * 3^3
\put {$\bullet$} at 25.5 -10
 \put {$\bullet$} at 25 -9
 \put {$\bullet$} at 26 -9
 \put {$\bullet$} at 26.5 -8
\plot 25.5 -10 25 -9 /
\plot 25.5 -10 26 -9 /
\plot 26 -9 26.5 -8 /

% 55 = 5*11
\put {$\bullet$} at 28 -10
 \put {$\bullet$} at 27.5 -9
 \put {$\bullet$} at 28.5 -9
\plot 28 -10 27.5 -9 /
\plot 28 -10 28.5 -9 /

% 56 = 2^3 * 7
\put {$\bullet$} at 30.5 -10
 \put {$\bullet$} at 30 -9
 \put {$\bullet$} at 29.5 -8
 \put {$\bullet$} at 31 -9
\plot 30.5 -10 30 -9 /
\plot 30 -9 29.5 -8 /
\plot 30.5 -10 31 -9 /

% 57 = 3*19
\put {$\bullet$} at 0.5 -15
 \put {$\bullet$} at 0 -14
 \put {$\bullet$} at 1 -14
\plot 0.5 -15 0 -14 /
\plot 0.5 -15 1 -14 /

% 58 = 2*29
\put {$\bullet$} at 2.5 -15
 \put {$\bullet$} at 2 -14
 \put {$\bullet$} at 3 -14
\plot 2.5 -15 2 -14 /
\plot 2.5 -15 3 -14 /

% 59
\put {$\bullet$} at 4 -15
 \put {$\bullet$} at 4 -14
\plot 4 -15 4 -14 /

% 60 = 2^2 * 3 * 5
\put {$\bullet$} at 7 -15
\put {$\bullet$} at 6 -14
\put {$\bullet$} at 5 -13
\put {$\bullet$} at 7 -14
\put {$\bullet$} at 8 -14
\plot 7 -15 6 -14 /
\plot 6 -14 5 -13 /
\plot 7 -15 7 -14 /
\plot 7 -15 8 -14 /

% 61
\put {$\bullet$} at 9 -15
 \put {$\bullet$} at 9 -14
\plot 9 -15 9 -14 /

% 62 = 2*31
\put {$\bullet$} at 10.5 -15
 \put {$\bullet$} at 10 -14
 \put {$\bullet$} at 11 -14
\plot 10.5 -15 10 -14 /
\plot 10.5 -15 11 -14 /

% 63 = 3^2 * 7
\put {$\bullet$} at 13 -15
 \put {$\bullet$} at 12.5 -14
 \put {$\bullet$} at 12 -13
 \put {$\bullet$} at 13.5 -14
\plot 13 -15 12.5 -14 /
\plot 12.5 -14 12 -13 /
\plot 13 -15 13.5 -14 /

% 64 = 2^{2*3}
 \put {$\bullet$} at 15 -15
 \put {$\bullet$} at 15 -14
 \put {$\bullet$} at 14.5 -13
 \put {$\bullet$} at 15.5 -13
 \plot 15 -15 15 -14 /
\plot 15 -14 14.5 -13 /
\plot 15 -14 15.5 -13 /

% 65 = 5*13
\put {$\bullet$} at 17 -15
 \put {$\bullet$} at 16.5 -14
 \put {$\bullet$} at 17.5 -14
\plot 17 -15 16.5 -14 /
\plot 17 -15 17.5 -14 /

% 66 = 6*11
\put {$\bullet$} at 19 -15
 \put {$\bullet$} at 18.5 -14
 \put {$\bullet$} at 19.5 -14
\plot 19 -15 18.5 -14 /
\plot 19 -15 19.5 -14 /

% 67
\put {$\bullet$} at 20.5 -15
 \put {$\bullet$} at 20.5 -14
\plot 20.5 -15 20.5 -14 /

% 68 = 2^2 * 17
\put {$\bullet$} at 22.5 -15
 \put {$\bullet$} at 22 -14
 \put {$\bullet$} at 21.5 -13
 \put {$\bullet$} at 23 -14
\plot 22.5 -15 22 -14 /
\plot 22 -14 21.5 -13 /
\plot 22.5 -15 23 -14 /

% 69 = 3*23
\put {$\bullet$} at 24.5 -15
 \put {$\bullet$} at 24 -14
 \put {$\bullet$} at 25 -14
\plot 24.5 -15 24 -14 /
\plot 24.5 -15 25 -14 /

% 70 = 2*5*7
\put {$\bullet$} at 27 -15
\put {$\bullet$} at 26 -14
\put {$\bullet$} at 27 -14
\put {$\bullet$} at 28 -14
\plot 27 -15 26 -14 /
\plot 27 -15 27 -14 /
\plot 27 -15 28 -14 /

% 71
\put {$\bullet$} at 29 -15
 \put {$\bullet$} at 29 -14
\plot 29 -15 29 -14 /

%72 = 2^3 * 3^2
\put {$\bullet$} at 31 -15
 \put {$\bullet$} at 30.5 -14
  \put {$\bullet$} at 30 -13
 \put {$\bullet$} at 31.5 -14
 \put {$\bullet$} at 32 -13
\plot 31 -15 30.5 -14 /
\plot 30.5 -14 30 -13 /
\plot 31 -15 31.5 -14 /
\plot 31.5 -14 32 -13 /

% 73
\put {$\bullet$} at 0 -20
 \put {$\bullet$} at 0 -19
\plot 0 -20 0 -19 /

% 74 = 2*37
\put {$\bullet$} at 1.5 -20
 \put {$\bullet$} at 1 -19
 \put {$\bullet$} at 2 -19
\plot 1.5 -20 1 -19 /
\plot 1.5 -20 2 -19 /

% 75 = 3 * 5^2
\put {$\bullet$} at 3.5 -20
 \put {$\bullet$} at 3 -19
 \put {$\bullet$} at 4 -19
 \put {$\bullet$} at 4.5 -18
\plot 3.5 -20 3 -19 /
\plot 3.5 -20 4 -19 /
\plot 4 -19 4.5 -18 /

% 76 = 2^2 * 19
\put {$\bullet$} at 6.5 -20
 \put {$\bullet$} at 6 -19
 \put {$\bullet$} at 5.5 -18
 \put {$\bullet$} at 7 -19
\plot 6.5 -20 6 -19 /
\plot 6 -19 5.5 -18 /
\plot 6.5 -20 7 -19 /

% 77 = 7*11
\put {$\bullet$} at 8.5 -20
 \put {$\bullet$} at 8 -19
 \put {$\bullet$} at 9 -19
\plot 8.5 -20 8 -19 /
\plot 8.5 -20 9 -19 /

% 78 = 2*3*13
\put {$\bullet$} at 11 -20
\put {$\bullet$} at 10 -19
\put {$\bullet$} at 11 -19
\put {$\bullet$} at 12 -19
\plot 11 -20 10 -19 /
\plot 11 -20 11 -19 /
\plot 11 -20 12 -19 /

% 79
\put {$\bullet$} at 13 -20
 \put {$\bullet$} at 13 -19
\plot 13 -20 13 -19 /

% 80 = 2^{2^2} * 5
\put {$\bullet$} at 15.5 -20
\put {$\bullet$} at 15 -19
\put {$\bullet$} at 14.5 -18
\put {$\bullet$} at 14 -17
\put {$\bullet$} at 16 -19
\plot 15.5 -20 15 -19 /
\plot 15 -19 14.5 -18 /
\plot 14.5 -18 14 -17 /
\plot 15.5 -20 16 -19 /

% 81 = 3^{2^2}
\put {$\bullet$} at 17 -20
\put {$\bullet$} at 17 -19
\put {$\bullet$} at 17 -18
\put {$\bullet$} at 17 -17
\plot 17 -20 17 -19 /
\plot 17 -19 17 -18 /
\plot 17 -18 17 -17 /

% 82 = 2*41
\put {$\bullet$} at 18.5 -20
 \put {$\bullet$} at 18 -19
 \put {$\bullet$} at 19 -19
\plot 18.5 -20 18 -19 /
\plot 18.5 -20 19 -19 /

% 83
\put {$\bullet$} at 20 -20
 \put {$\bullet$} at 20 -19
\plot 20 -20 20 -19 /

% 84 = 2^2 * 3 * 7
\put {$\bullet$} at 23 -20
\put {$\bullet$} at 22 -19
\put {$\bullet$} at 21 -18
\put {$\bullet$} at 23 -19
\put {$\bullet$} at 24 -19
\plot 23 -20 22 -19 /
\plot 22 -19 21 -18 /
\plot 23 -20 23 -19 /
\plot 23 -20 24 -19 /

% 85 = 5*17
\put {$\bullet$} at 25.5 -20
 \put {$\bullet$} at 25 -19
 \put {$\bullet$} at 26 -19
\plot 25.5 -20 25 -19 /
\plot 25.5 -20 26 -19 /

% 86 = 2*43
\put {$\bullet$} at 27.5 -20
 \put {$\bullet$} at 27 -19
 \put {$\bullet$} at 28 -19
\plot 27.5 -20 27 -19 /
\plot 27.5 -20 28 -19 /

% 87 = 3*29
\put {$\bullet$} at 29.5 -20
 \put {$\bullet$} at 29 -19
 \put {$\bullet$} at 30 -19
\plot 29.5 -20 29 -19 /
\plot 29.5 -20 30 -19 /

% 88 = 2^3 * 11
\put {$\bullet$} at 1 -25
 \put {$\bullet$} at 0.5 -24
 \put {$\bullet$} at 0 -23
 \put {$\bullet$} at 1.5 -24
\plot 1 -25 0.5 -24 /
\plot 0.5 -24 0 -23 /
\plot 1 -25 1.5 -24 /

% 89
\put {$\bullet$} at 2.5 -25
 \put {$\bullet$} at 2.5 -24
\plot 2.5 -25 2.5 -24 /

% 90 = 2 * 3^2 * 5
\put {$\bullet$} at 4.5 -25
\put {$\bullet$} at 3.5 -24
\put {$\bullet$} at 4.5 -24
\put {$\bullet$} at 5.5 -24
\put {$\bullet$} at 4.5 -23
\plot 4.5 -25 3.5 -24 /
\plot 4.5 -25 4.5 -24 /
\plot 4.5 -25 5.5 -24 /
\plot 4.5 -24 4.5 -23 /

% 91 = 7*13
\put {$\bullet$} at 7 -25
 \put {$\bullet$} at 6.5 -24
 \put {$\bullet$} at 7.5 -24
\plot 7 -25 6.5 -24 /
\plot 7 -25 7.5 -24 /

% 92 = 2^2 * 23
\put {$\bullet$} at 9.5 -25
 \put {$\bullet$} at 9 -24
 \put {$\bullet$} at 8.5 -23
 \put {$\bullet$} at 10 -24
\plot 9.5 -25 9 -24 /
\plot 9 -24 8.5 -23 /
\plot 9.5 -25 10 -24 /

% 93 = 3 * 31
\put {$\bullet$} at 11.5 -25
 \put {$\bullet$} at 11 -24
 \put {$\bullet$} at 12 -24
\plot 11.5 -25 11 -24 /
\plot 11.5 -25 12 -24 /

% 94 = 2 * 47
\put {$\bullet$} at 13.5 -25
 \put {$\bullet$} at 13 -24
 \put {$\bullet$} at 14 -24
\plot 13.5 -25 13 -24 /
\plot 13.5 -25 14 -24 /

% 95 = 5 * 19
\put {$\bullet$} at 15.5 -25
 \put {$\bullet$} at 15 -24
 \put {$\bullet$} at 16 -24
\plot 15.5 -25 15 -24 /
\plot 15.5 -25 16 -24 /

% 96 = 2^5 * 3
\put {$\bullet$} at 18 -25
 \put {$\bullet$} at 17.5 -24
 \put {$\bullet$} at 17 -23
 \put {$\bullet$} at 18.5 -24
\plot 18 -25 17.5 -24 /
\plot 17.5 -24 17 -23 /
\plot 18 -25 18.5 -24 /

% 97
\put {$\bullet$} at 19.5 -25
 \put {$\bullet$} at 19.5 -24
\plot 19.5 -25 19.5 -24 /

% 98 = 2 * 7^2
\put {$\bullet$} at 21 -25
 \put {$\bullet$} at 20.5 -24
 \put {$\bullet$} at 21.5 -24
 \put {$\bullet$} at 22 -23
\plot 21 -25 20.5 -24 /
\plot 21 -25 21.5 -24 /
\plot 21.5 -24 22 -23 /

% 99 = 3^2 * 11
\put {$\bullet$} at 24 -25
 \put {$\bullet$} at 23.5 -24
 \put {$\bullet$} at 23 -23
 \put {$\bullet$} at 24.5 -24
\plot 24 -25 23.5 -24 /
\plot 23.5 -24 23 -23 /
\plot 24 -25 24.5 -24 /

% 100 = 2^2 * 5^2
\put {$\bullet$} at 26.5 -25
 \put {$\bullet$} at 26 -24
  \put {$\bullet$} at 25.5 -23
 \put {$\bullet$} at 27 -24
 \put {$\bullet$} at 27.5 -23
\plot 26.5 -25 26 -24 /
\plot 26 -24 25.5 -23 /
\plot 26.5 -25 27 -24 /
\plot 27 -24 27.5 -23 /

% 101
\put {$\bullet$} at 28.5 -25
 \put {$\bullet$} at 28.5 -24
\plot 28.5 -25 28.5 -24 /

% 102 = 2*3*17
\put {$\bullet$} at 30.5 -25
\put {$\bullet$} at 29.5 -24
\put {$\bullet$} at 30.5 -24
\put {$\bullet$} at 31.5 -24
\plot 30.5 -25 29.5 -24 /
\plot 30.5 -25 30.5 -24 /
\plot 30.5 -25 31.5 -24 /

% \arrow <0.215cm> [0.2,0.6] from 0.5 1.5 to 3.5 2.5 

\endpicture \]

 \caption{${\cal T}({\mathbb N})$ in lexicographic order up to $t(102)$}
 \label{viale alberato}
\end{figure}

Planar rooted trees, or rooted trees, can also be identified by using the parentheses syntax (see \cite{HPT,DBM} and \cite{CCF} for a polynomial classification). Examples of the map $t$ are: 
\begin{equation} 
t(1)=() \; ,
\end{equation}

\begin{equation} 
t(2)=t(3)=t(p)=(()) \; , \quad \text{for all primes} \; p
\end{equation}

\begin{equation} 
t(12)=((())()) \; ,
\end{equation}

\begin{equation} 
t(18)=t(486)=(()(())) \; ,
\end{equation}

and, from the example $(\ref{154})$,

\begin{equation} 
t(1549681956)=((())((()()))) \; .
\end{equation}
More examples of the map $t(n)$, with both the graph representation and the parenthesis syntax, can be found in Figures $\ref{tiofn}$ and $\ref{moretiofn}$. In Figure $\ref{viale alberato}$ ${\cal T}({\mathbb N})$ is shown up to $t(102)$.

By the typical property of $t^\#$ we have, for instance, $t^\#(12)=t^\#(18)$ or $t^\#(48)=t^\#(162)$. We observe that by decorating unlabeled trees with suitable primes one sees that any finite planar (or general) rooted tree may appear, i.e. both the maps $t$ and $t^\#$ are surjective.

Listing the integers where new trees arise in ${\cal T}$ we find the sequence

\begin{equation}\label{tau}
1, 2, 4, 6, 12, 16, 18, 30, 36, 48, \ldots
\end{equation}

as well as the corresponding subsequence for $[{\cal T}]$
\begin{equation}\label{taud}
1, 2, 4, 6, 12, 16, 30, 36, 48, \ldots
\end{equation}
where we notice e.g. the disappearence of the integer $18$ that has, in $[{\cal T}]$, the same structure of $12$.
The two sequences are identified as {\it A284456} and {\it A279686} in \cite{OEIS}.

Our main object of interest will be the sequence of trees corresponding to the natural numbers:
\begin{equation}
{\cal T}({\mathbb N}) = \{t(n)\}_{n=1}^{\infty}\, .
\end{equation}
We will refer to it as the ${\mathbb N}$atural Avenue, since typically an avenue is a road with trees on the side at regular intervals. Similarly $[{\cal T}]({\mathbb N})$ will be the corresponding sequence of trees generated from the map $t^\#$.

\medskip
{\it 
We are ready now to state the general perspective that has inspired our interest and led us to write this note. We have obtained ${\cal T}({\mathbb N})$ from ${\mathbb N}$ by erasing the primes labeling the tree corresponding to each integer. In doing so we have lost some information. Is it possible to recover that information, i.e. to invert the map $t$, for instance by means of local observations? In other words, by observing finite windows of ${\cal T}({\mathbb N})$ can we discover where we are located? Moreover, concerning the global structure of ${\cal T}({\mathbb N})$, how are the different trees, or groups of them, distributed? What are their specific densities and fluctuations?} 
\medskip

In the following section we start the investigation of the local and the global structures of ${\cal T}({\mathbb N})$ and introduce some useful concepts.

\section{The structure of \texorpdfstring{${\cal T}({\mathbb N})$}{}}

Here we consider (finite) windows of ${\mathbb N}$ and we denote them as $[n,n+1,n+2,\cdots, n+l]$. The elements of ${\cal T}$, the planar rooted trees, can be used as letters of an alphabet to form (finite) {\it words}. For instance, $[(()),(())]$ and $[(()),(()),(())]$ are two words. While the first can be found 
in ${\cal T}({\mathbb N})$ only once, in correspondence of the window $[2,3]$, the second cannot be found at all because there are no three prime numbers of the form $p,p+1,p+2$. Words can also include, in their interior (i.e. not at the boundary), a special letter, the asterisk $*$, denoting the fact that any possible tree can appear in such a position. A similar word is, for instance, $[(()),*,*,*,(())]$ that appears in correspondence to the windows $[3,4,5,6,7]$, $[7,8,9,10,11]$, $[13,14,15,16,17]$, etc. 

When a word can be found in ${\cal T}({\mathbb N})$, like in the previous example, it is called a {\it natural} word or {\it observable} for short. On the other hand, when a word cannot be found in ${\cal T}({\mathbb N})$ it is called {\it unnatural}.

Observables can be suitably identified by one of their corresponding windows in ${\mathbb N}$. For example, the observable $O(3,5)$ represents the natural word $[(()),*,(())]$, namely the trees associated to two prime numbers and an arbitrary tree in between; 
of course, it corresponds precisely to 
a couple of {\it twin prime numbers}. A further observable is $O(3,4,5)$ i.e. the sequence [(()),((())),(())]. Note that $O(3,5)$ and $O(3,4,5)$ are different observables, the first including the asterisk in the middle while the second is a word without asterisks. The second one corresponds to a couple of twin prime numbers separated by a prime raised to a prime power.

It is straightforward to obtain some very basic structural properties of ${\cal T}({\mathbb N})$ starting from elementary facts about the decomposition of an integer in primes. Other properties are more challenging and can been proved with advanced mathematical methods. In some other case, these properties are so difficult that a proof is still lacking (some of them are well known classical conjectures). To give a flavour of all that, we investigate the structure of some observables appearing nearby the origin. Denoting by $|O|$ the number of occurrences of the observable $O$ in the ${\mathbb N}$atural Avenue we have the following lemmata and open questions:
\begin{itemize}
\item[-] $|O(1)|=1$, 
\item[-] $|O(2)| = \infty$, i.e. the observable $[(())]$ appears infinitely many times as a consequence of the Euclid theorem on the existence of infinitely many primes. 
\item[-] $|O(n)|=\infty$ for $n\geq 2$, by suitably decorating the tree corresponding to $n$ with primes.
\item[-] $|O(2,3)| 
= 1$, as $(2,3)$ is the only contiguous couple of prime numbers.
\item[-] $|O(4,5)| = 1$, by using the divisibility properties by 2 and 3. Indeed, we have to solve the equation $2^p + 1 = q$ with $p$ and $q$ prime numbers. Now, if $p >2$ we have that $(3-1)^p + 1$ is always divisible by 3.
\item[-] $|O(3,5)| 
= \infty ?$, This amounts to the {\it twin prime conjecture}, a very well-known open problem (although important progress has been achieved recently \cite{YZ}, see also \cite{Ma} for a recent account).
\item[-] $|O(3,5,7)| = 1$, observing that if $p, p+2,p+4$ are primes then the corresponding congruence classes {\rm mod} 3 exhaust all possible cases, 
i.e.  $ \{p,p+2,p+4\} = \{0,1,2\} \; ({\rm mod} \ 3)$ (up to reordering). Thus one of these three primes, being divisible by 3, must actually be equal to 3. The only possible case is that indeed $p=3$.

\item[-] $|O(3,4)| = \infty ?$ Indeed this observable selects precisely the {\it Mersenne numbers}, i.e. primes of the form $2^p-1$ for some prime $p$ (these are the first entries of each pair). Again, the cardinality of the set of Mersenne numbers is another long-standing open problem.

\item[-] $|O(5,6)| = \infty ?$ This observable selects primes $p$ such that $(p+1)/2$ is also prime or, equivalently, primes $q$ such that $2q - 1$ is also prime (see A005383 in OEIS). It is unknown to us whether there are infinitely many such primes.
\item[-] $|O(6,7)| = \infty?$ This observable selects the {\it Sophie Germain primes}, i.e. primes $q$ such that $2q + 1$ is also prime (see A005385 in OEIS).
It is unknown whether there are infinitely many such primes.

\item[-] $|O(5,6,7)| = 1$ as easily checked, for we have triples of consecutive integers of the form $(p,qr,s)$, with $p,q,r,s$ primes and $q\neq r$, which can only happen when $q=2$ and $r=3$.
\item[-] $|O(8,9)|=1$. This amounts to finding the solutions of 
$2^p \pm 1=q^r$ with $p,q,r$ primes and was proved in \cite{PM} in relation to the Catalan conjecture. 
\item[-] $|O(11,12)| = \infty?$ This amounts to finding the solutions of $p+1=2^q r$ with $p,q,r$ primes with $r\neq 2$. The sequence of the $p$ starts with $11,19,23,43,67,103,151,$ $163,211,223,283,\ldots$
 which is not listed in OEIS. It is unknown to us whether there are infinitely many such primes.
\item[-] $|O(12,13)| = \infty?$ This amounts to finding the solutions of $p-1=2^q r$ with $p,q,r$ primes with $r\neq 2$. The sequence of the $p$ starts with $13,29,41,53,89,97,137,$ $149,173,233,269,\ldots$
 which is not listed in OEIS. It is unknown to us whether there are infinitely many such primes.
\item[-] $|O(11,12,13)| = |O(17,18,19)| = 1$.
Indeed, any other occurrence of the first observable should then be of the form $O(q,2^p 3, q+2)$, with $q,q+2,p$ prime.
However, for all primes $p \geq 3$, we claim that at least one among the two integers $2^p 3 \pm 1$ is always divisible by 5 (and actually only one of them). To see this, write $2^p 3 \pm 1 = 24 \cdot 2^{p-3} \pm 1 = (5^2-1) 2^{p-3} \pm 1$. It follows that one among $2^p 3 \pm 1$ is divisible by 5 if and only if one among $-2^{p-3} \pm 1$ is divisible by 5 or, equivalently, one among $2^{p-3} \pm 1$ is divisible by 5. Also notice that the statement we want to prove is certainly true for $p = 3$, so in the sequel we can assume that $p \geq 5$ without loss of generality. Consider the five consecutive integers
$$2^{p-3} - 2, \ 2^{p-3} - 1, \ 2^{p-3}, \ 2^{p-3} + 1, \ 2^{p-3} + 2 \ ,$$
then one and only one of them is divisible by 5. But clearly $2^{p-3}$ is not divisible by $5$, so our claim will be proved if we show that both $2^{p-3} \pm 2$ are not divisible by 5.
Write $2^{p-3} \pm 2 = 2(2^{p-4} \pm 1)$. Thus, it suffices to show that both $2^{p-4} \pm 1$ are not divisible by 5.
Since $p$ is odd, we can write $2^{p-4} \pm 1 = 2 (2^{p-5}) \pm 1 = 2 \cdot 4^{(p-5)/2} \pm 1 = 2 (5-1)^{(p-5)/2} \pm 1$, which will be divisible by 5 if and only if $2 (-1)^{(p-5)/2} \pm 1$ is divisible by 5. But the latter condition is clearly 
not satisfied, and the result for the first observable is proved. As for the second, for $p \geq 3$ prime, $2^p \cdot 3 \pm 1 = (5-3)^p (5-2) \pm 1$ is divisible by 5 if and only if $-2 (-3)^p \pm 1 = 2 \cdot 3^p \pm 1$ is divisible by 5, which concludes the proof. Notice that the two integers $12$ and $18$ correspond to different planar trees belonging to the same non-planar class of equivalence.
\item[-] $|O(3,5)|+|O(2,4)|+|O(13,15)|=\infty$. This is a consequence of the result proved in \cite{CJ} stating that there are infinitely many primes $p$ (Chen primes) such that $p+2$ is either a prime or a product of two primes not necessarily distinct. In our context this tells us that at least one of the three involved observables appears infinitely many times in ${\cal T}({\mathbb N})$.
\item[-] 
Similarly, another example of an observable appearing infinitely many times that is not a single tree arises from 
Zhang's theorem \cite{YZ}, 
namely the existence of infinitely many primes differing at most by a given (possibly large, but finite) number $C$:
choosing finitely many pairs of primes $(p_i,q_i)$ such that $q_i - p_i = 2 i$, $i = 1,\ldots, [C/2]+1$,
we get
$$\sum_{i=1}^{[C/2] + 1} |O(p_i,q_i)| = +\infty \ $$
(in the sum we can omit those $i$'s for which the pair $(p_i,q_i)$ as above does not exist).
\end{itemize}

\medskip

In general, an observable can appear a finite or an infinite number of times in ${\cal T}({\mathbb N})$. If it is unique, i.e. it appears only once, it provides the complete information of its location. For instance, if we observe $O(4,5)$ or $O(8,9)$ we know for sure where we are placed in ${\mathbb N}$, while the observation of a single tree does not pinpoint our position. If an observable appears only once then every other observable that includes it (as a subword) appears only once as well. We will call {\it milestones} those observables that appear only once and are irreducible, i.e. do not properly contain any other observable that appears only once. For instance $O(4,5,6)$ appears only once but it is not a milestone because it properly contains $O(4,5)$. The latter observable is a milestone because $O(4)$ and $O(5)$ appear infinitely many times.

\subsection{There are infinitely many milestones}
Given $n \in {\mathbb N}$, let us define $\kappa_+(n)$ as the least non-negative integer such that $O(n,n+1,\ldots,n+\kappa_+(n))$ is unique. 
Although it may not be obvious that such an integer always exists, we show below that this is indeed the case. 
We also believe that $\kappa_+$ is an unbounded function.
It readily follows from the above discussion that $\kappa_+(1) = 0$, $\kappa_+(2) = 1$, $\kappa_+(3) = 2$ (since $O(4,5)$ and thus $O(3,4,5)$ are unique)
and $\kappa_+(4) = 1$. 
We also see that $\kappa_+(5) = 2$ as well (note that $O(5,6) = O(13,14)$).

Defining $\kappa_-$ in a similar manner but looking at backward windows, we clearly have $\kappa_-(n) \leq n-2$, for $n \geq 3$. Moreover $\kappa_-(1) = 0$, $\kappa_-(2) = 1$, $\kappa_-(3) = 1$. Therefore, $\kappa_+$ and $\kappa_-$ are different and actually one is not even bounded by the other in a strict sense because $\kappa_-(4) = 2$, $\kappa_-(5) = 1$.

We now address the problems of the finiteness of $\kappa_+(n)$ and of the cardinality of milestones.
To prove that $\kappa_+(n)$ is finite for every $n$ we use a strategy that generalizes the proof of the uniqueness of the observable $O(3,5,7)$.

\begin{theorem}
For each natural $n$ there exists $\kappa_+(n) \in {\mathbb N}$ 
with the properties
that the observable $O(n,\ldots,n+\kappa_+(n))$ is unique and moreover $\kappa_+(n)$ is the least such number.    
\end{theorem}

\begin{proof} 
Let $n$ be a natural number, and let $p = p_0$ be the first prime larger or equal to $n$. Thanks to Dirichlet theorem on primes in arithmetic progressions \cite{AT},
can find a prime $p_1$ such that $p_1 - p_0 > p_0$ and $p_1 \cong 1 \; ({\rm mod} \ p)$. 
Likewise, for every $ 2 \leq i \leq p-1$ can also find $p_{i} > p_{i-1}$ such that $p_i \cong i \; ({\rm mod} \ p)$.
If $O(p_0,*,\ldots,*,p_1,*,\ldots,*,p_{p-1})$ were not unique, say 
$O(p_0,*,\ldots,*,p_1,*,\ldots,*,p_{p-1}) = O(p'_0,*,\ldots,*,p'_1,*,\ldots,*,p'_{p-1})$ 
then one of the $p'_i$, $0 \leq i \leq p-1$ should be (divisible and thus) equal to $p$, and it is easy to see that the only possible such occurrence is indeed $p'_0 = p$.
Thus, $O(p_0,*,\ldots,*,p_1,*,\ldots,*,p_{p-1})$ is unique, and hence $O(n,n+1,\ldots,p_{p-1})$ is unique as well. The conclusion follows at once.
\end{proof}
As a consequence, we immediately obtain the following important result.
\begin{corollary}
There are infinitely many milestones.          
\end{corollary}

\medskip
So far we mostly focused on ``local'' observables.
However, the global features of ${\cal T}({\mathbb N})$ are also worth of investigation. 
For instance it is not difficult to see that the ``infinite observable'' $O_+[n]:= O(n,n+1,n+2,\ldots)$ completely determines $n$, i.e. if two natural numbers $n$ and $m$ are such that $O_+[n] = O_+[m]$ then $n = m$. Of course this is an obvious consequence of the previous theorem but an independent elementary proof by contradiction goes as follows. Suppose, without loss of generality, that $m \geq n$ and set $r:= m-n$. Then it readily follows from the above assumption that eventually $p + r$ is prime whenever $p$ is prime (this is a consequence of the fact that every tree appearing in the $h$-th slot of $O_+(n)$ will also appear in the corresponding slot of $O_+(m)$).
Hence, $r$ must be equal to 0, for $r >0$ would contradict the well-known fact that $\limsup_k (p_{k+1} - p_k) = + \infty$ (where $p_k$ is the $k$-th prime number). 

Other properties addressing global aspects of ${\cal T}({\mathbb N})$ will appear in the next section.\\

\medskip
For each observable $O$, we define $O^t$, the {\it transposed} of $O$, as the word obtained by reversing the sequence of trees $O$. For example $O^t(5,6)=O(6,7)$. It is not always true that the transposed of an observable is itself an observable, this happens for instance with $O(1,2) = [(),(())]$. As one may expect the transposition changes in general the cardinality of the occurrences, for instance $O(3,4)$ appears many times while $O(4,5)=O^t(3,4)$ is unique.
\medskip

\section{Some perspectives}\label{persp}

Below we will lay down some basic directions about the way we plan to
examine ${\cal T}({\mathbb N})$, and state some problems/conjectures. We start with the elementary observation that for each observable that appears infinitely many times there is a corresponding problem of determining its distribution (density and fluctuations).

In particular this applies to the observables made of a single tree, since every tree clearly appears infinitely many times. This problem in its very basic instance, the first nontrivial tree, 
is the celebrated problem of the prime number distribution, the most fundamental and fascinating topic in analytic number theory. Indicating by $\pi(x)$ the number of prime numbers less or equal than $x$, Jacques Hadamard and 
Charles Jean de la Vall\'ee Poussin proved in 1896 that 
$$
\lim_{x\to \infty}\frac{\pi(x)}{\frac{x}{\log x}}=1 \; ,
$$
a result referred to as the {\it prime number theorem}. The, still unproven, Riemann hypothesis would imply an accurate estimate of the fluctuations around that density, namely for large $x$
$$
\pi(x)={\rm li}(x)+O(\sqrt{x}\log x)\; .
$$
The prime number theorem suggests that the spacing between consecutive occurrences of primes is unbounded. Actually this can indeed be proven by very elementary methods (see for instance \cite{TTB1}, line 11 for this fact, and much more).

From our perspective we have an infinite family of density and fluctuation problems and we expect that some relations should appear between the densities of different trees (see \cite{Na}). The list includes some well-know classical problems 
that have been solved more than a century ago. They deal e.g. with the distribution of products of two distinct primes (``square-free semiprimes'') or, more generally, of square-free integers. These problems correspond in our context to the distribution of the {\it bushes} of the form $[(())(())]$ or it's n-fold generalization $[(())\cdots(())]$ respectively.
It turns out that the number of products of distinct primes $p$ and $q$ with $p<q$, $pq \leq x$ goes, for large $x$, like $x \log\log x / \log x$ \cite{La1} (see also \cite{La2}).
The asymptotic density of all the square free integers is $6/\pi^2 \approx 0.6079$ \cite{Ge}.
See also a recent work on some asymptotic formulas related to the height of the trees in \cite{VVI2}.

\medskip
Summarizing the above discussion the natural questions that arise from our perspective can be grouped into three areas. {\it Multiplicities}: to our knowledge the observables in ${\cal T}({\mathbb N})$ take only multiplicities 1 or infinity. Are there also observables with finite multiplicity larger than 1? We remind that instead in $[{\cal T}]({\mathbb N})$ we 
found an observable of cardinality 2. 
{\it Distributions}: compute (by rigorous arguments, but also empirically guess using some high-performance machine) the distribution (densities and fluctuations) of observables with infinite multiplicity. Find possible relations between such densities. {\it Spacing}: considering an observable with infinite multiplicity, are the spacings between consecutive appearances unbounded like it happens for the primes?

We have also shown that some of the questions that naturally emerge from our considerations are related to classical results or conjectures in number theory. But our perspective poses many more questions, and we hope that they will be useful to shed further light on the overall subject. We plan to return on some of these topics in the future.

\section*{Acknowledgements}
v1: This work was basically completed, in the form of handwritten notes that we still have, already in 1992 (with the only exception 
of the recent bibliography) while the two authors were working on their PhD theses in operator algebras and statistical mechanics, respectively. 
The discussions around the themes presented here had at that time, and still have now, no other purpose than to set our curiosity free and the precious side effect to distract our minds from the urgent academic duties. It's a pleasure to thank the Physics Department of the University of Rome La Sapienza for the generous hospitality to let us use, at that time, the seminar room {\it Enrico Persico}. \\
v2: We thank R. Schoof for a crucial suggestion that led us to the proof of the existence of infinitely many milestones.

\newpage

\begin{figure}
 \begin{subfigure}[b]{0.24\textwidth}
\centering
\begin{tikzpicture}
\tikzstyle{hollow node}=[circle,draw,inner sep=1.5]
\tikzstyle{solid node}=[circle,draw,inner sep=1.5,fill=black]
\node(0)[solid node]{};
%child[grow=90]{node[hollow node]{}};
\end{tikzpicture}
\caption{$t(1)=()$}
    \label{1}
    \end{subfigure}
    \hfill
     \begin{subfigure}[b]{0.24\textwidth}
     \centering
\begin{tikzpicture}
\tikzstyle{hollow node}=[circle,draw,inner sep=1.5]
\tikzstyle{solid node}=[circle,draw,inner sep=1.5,fill=black]
\node(0)[solid node]{}
child[grow=90]{node[hollow node]{}};
\end{tikzpicture}
\caption{$t(2)=(())$}
    \label{2}
     \end{subfigure}
     \hfill
     \begin{subfigure}[b]{0.24\textwidth}
     \centering
\begin{tikzpicture}
\tikzstyle{hollow node}=[circle,draw,inner sep=1.5]
\tikzstyle{solid node}=[circle,draw,inner sep=1.5,fill=black]
\node(0)[solid node]{}
child[grow=90]{node[hollow node]{}};
\end{tikzpicture}
\caption{$t(3)=(())$}
    \label{3}
     \end{subfigure}
    \hfill 
     \begin{subfigure}[b]{0.24\textwidth}
     \centering
\begin{tikzpicture}
\tikzstyle{hollow node}=[circle,draw,inner sep=1.5]
\tikzstyle{solid node}=[circle,draw,inner sep=1.5,fill=black]
\node(0)[solid node]{}
child[grow=90]{node[hollow node]{}child[grow=90]{node[hollow node]{}}};
\end{tikzpicture}
\caption{$t(4)=((()))$}
    \label{4}
     \end{subfigure}
     %\vspace{0.5cm}
     \begin{subfigure}[b]{0.24\textwidth}
     \centering
\begin{tikzpicture}
\tikzstyle{hollow node}=[circle,draw,inner sep=1.5]
\tikzstyle{solid node}=[circle,draw,inner sep=1.5,fill=black]
\node(0)[solid node]{}
child[grow=90]{node[hollow node]{}};
\end{tikzpicture}
\caption{$t(5)=(())$}
    \label{5}
     \end{subfigure}
     \hfill
     \begin{subfigure}[b]{0.24\textwidth}
     \centering
\begin{tikzpicture}
\tikzstyle{hollow node}=[circle,draw,inner sep=1.5]
\tikzstyle{solid node}=[circle,draw,inner sep=1.5,fill=black]
\node(0)[solid node]{}
child[grow=60]{node[hollow node]{}}
child[grow=120]{node[hollow node]{}};
\end{tikzpicture}
\caption{$t(6)=(()())$}
    \label{6}
     \end{subfigure}
     \hfill
     \begin{subfigure}[b]{0.24\textwidth}
     \centering
\begin{tikzpicture}
\tikzstyle{hollow node}=[circle,draw,inner sep=1.5]
\tikzstyle{solid node}=[circle,draw,inner sep=1.5,fill=black]
\node(0)[solid node]{}
child[grow=90]{node[hollow node]{}};
\end{tikzpicture}
\caption{$t(7)=(())$}
    \label{7}
     \end{subfigure}
     \hfill
     \begin{subfigure}[b]{0.24\textwidth}
     \centering
     \vspace{0.5cm}
\begin{tikzpicture}
\tikzstyle{hollow node}=[circle,draw,inner sep=1.5]
\tikzstyle{solid node}=[circle,draw,inner sep=1.5,fill=black]
\node(0)[solid node]{}
child[grow=90]{node[hollow node]{}child[grow=90]{node[hollow node]{}}};
\end{tikzpicture}
\caption{$t(8)=((()))$}
    \label{8}
     \end{subfigure}
     
     \begin{subfigure}[b]{0.24\textwidth}
     \centering
     \vspace{0.5cm}
\begin{tikzpicture}
\tikzstyle{hollow node}=[circle,draw,inner sep=1.5]
\tikzstyle{solid node}=[circle,draw,inner sep=1.5,fill=black]
\node(0)[solid node]{}
child[grow=90]{node[hollow node]{}child[grow=90]{node[hollow node]{}}};
\end{tikzpicture}
\caption{$t(9)=((()))$}
    \label{9}
     \end{subfigure}
     \hfill
     \begin{subfigure}[b]{0.24\textwidth}
     \centering
\begin{tikzpicture}
\tikzstyle{hollow node}=[circle,draw,inner sep=1.5]
\tikzstyle{solid node}=[circle,draw,inner sep=1.5,fill=black]
\node(0)[solid node]{}
child[grow=60]{node[hollow node]{}}
child[grow=120]{node[hollow node]{}};
\end{tikzpicture}
\caption{$t(10)=(()())$}
    \label{10}
     \end{subfigure}
     \hfill
     \begin{subfigure}[b]{0.24\textwidth}
     \centering
\begin{tikzpicture}
\tikzstyle{hollow node}=[circle,draw,inner sep=1.5]
\tikzstyle{solid node}=[circle,draw,inner sep=1.5,fill=black]
\node(0)[solid node]{}
child[grow=90]{node[hollow node]{}};
\end{tikzpicture}
\caption{$t(11)=(())$}
    \label{11}
     \end{subfigure}
     \hfill
     \begin{subfigure}[b]{0.24\textwidth}
     \centering
\begin{tikzpicture}
\tikzstyle{hollow node}=[circle,draw,inner sep=1.5]
\tikzstyle{solid node}=[circle,draw,inner sep=1.5,fill=black]
\node(0)[solid node]{}
child[grow=60]{node[hollow node]{}}
child[grow=120]{node[hollow node]{}child[grow=120]{node[hollow node]{}}};
\end{tikzpicture}
\caption{$t(12)=((())())$}
    \label{12}
     \end{subfigure}
     \caption{The map $t(n)$ for $1\le n \le 12$.}
     \label{tiofn}
  \end{figure}
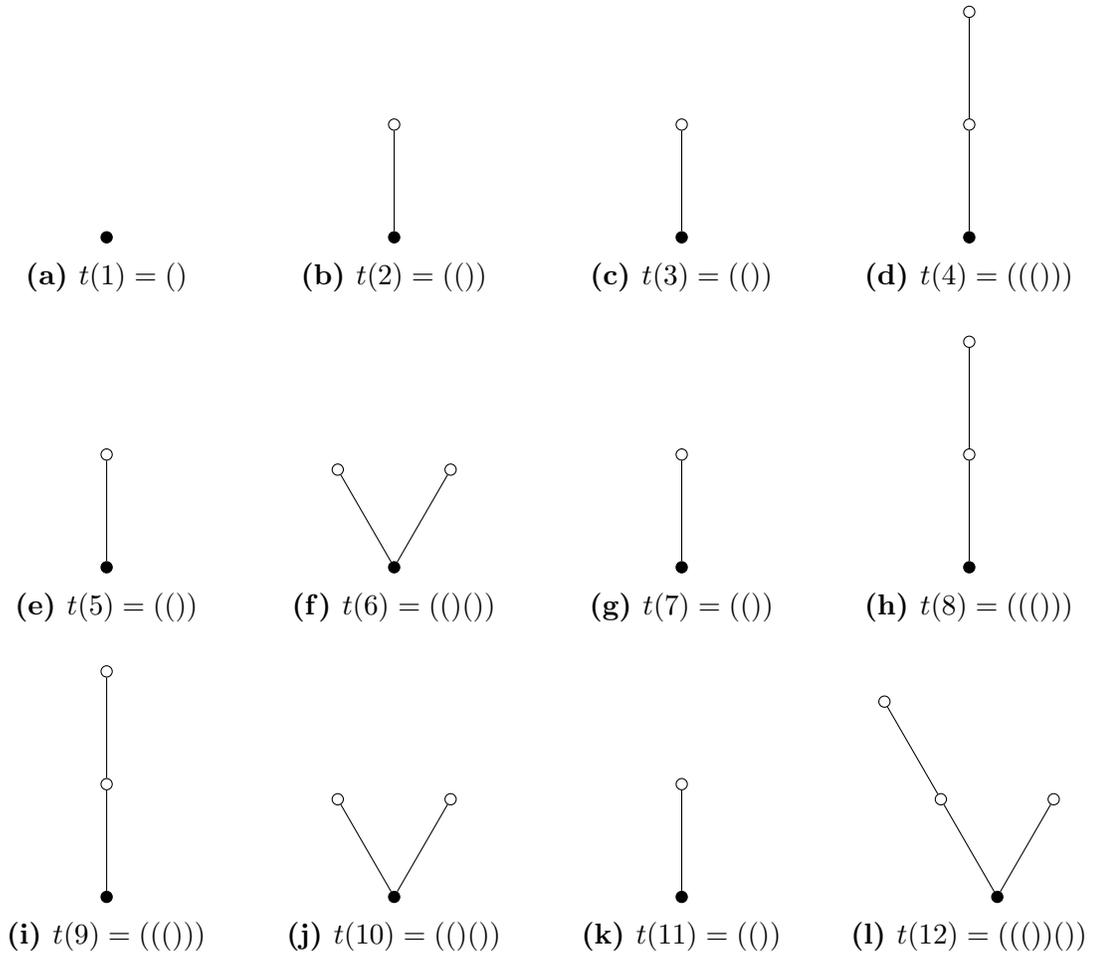 

%\end{document}

\end{document}